\documentclass[12pt]{article}
\usepackage{microtype} \DisableLigatures{encoding = *, family = * }
\usepackage{subeqn,multirow}
\usepackage[pdftex]{graphicx}
\usepackage{amsfonts,amssymb}
\usepackage{natbib}
\bibpunct{(}{)}{;}{a}{,}{,}
\usepackage[bookmarksnumbered=true, pdfauthor={Wen-Long Jin}]{hyperref}

\newcommand{\commentout}[1]{}

\newcommand{\ba}{\begin{array}}
        \newcommand{\ea}{\end{array}}
\newcommand{\bc}{\begin{center}}
        \newcommand{\ec}{\end{center}}
\newcommand{\bdm}{\begin{displaymath}}
        \newcommand{\edm}{\end{displaymath}}
\newcommand{\bds} {\begin{description}}
        \newcommand{\eds} {\end{description}}
\newcommand{\ben}{\begin{enumerate}}
        \newcommand{\een}{\end{enumerate}}
\newcommand{\beq}{\begin{equation}}
        \newcommand{\eeq}{\end{equation}}
\newcommand{\bfg} {\begin{figure}}
        \newcommand{\efg} {\end{figure}}
\newcommand{\bi} {\begin {itemize}}
        \newcommand{\ei} {\end {itemize}}
\newcommand{\bpp}{\begin{pspicture}}
        \newcommand{\epp}{\end{pspicture}}
\newcommand{\bqn}{\begin{eqnarray}} 
        \newcommand{\eqn}{\end{eqnarray}}
\newcommand{\bqs}{\begin{eqnarray*}}
        \newcommand{\eqs}{\end{eqnarray*}}
\newcommand{\bsq}{\begin{subequations}}
        \newcommand{\esq}{\end{subequations}}
\newcommand{\bsl} {\begin{slide}[8.8in,6.7in]}
        \newcommand{\esl} {\end{slide}}
\newcommand{\bss} {\begin{slide*}[9.3in,6.7in]}
        \newcommand{\ess} {\end{slide*}}
\newcommand{\btb} {\begin {table}[h]}
        \newcommand{\etb} {\end {table}}
\newcommand{\bvb}{\begin{verbatim}}
        \newcommand{\evb}{\end{verbatim}}

\newcommand{\m}{\mbox}
\newcommand {\der}[2] {{\frac {\m {d} {#1}} {\m{d} {#2}}}}

\newcommand{\cas}[1]{{{\left \{ \ba #1 \ea \right. }}}

\newcommand{\reff}[1] {{{Figure \ref {#1}}}}
\newcommand{\refe}[1] {{{(\ref {#1})}}}
\newcommand{\reft}[1] {{{\textbf{Table} \ref {#1}}}}

\newtheorem{theorem}{Theorem}[section]
\newtheorem{definition}[theorem]{Definition}
\newtheorem{lemma}[theorem]{Lemma}
\newtheorem{corollary}[theorem]{Corollary}

\usepackage{slashbox}

\def\La	{{\Lambda}}
\def\pmb#1{\setbox0=\hbox{$#1$}%
   \kern-.025em\copy0\kern-\wd0
   \kern.05em\copy0\kern-\wd0
   \kern-.025em\raise.0433em\box0 }

\def\eop{{\hfill $\blacksquare$}}
\def\r{{\rho}}

\def\dt     {{\Delta t}}

\def\eop{{\hfill $\blacksquare$}}

\def\la {{{\lambda}}}

\usepackage{cleveref}

\begin {document}
\title{Point queue models: a unified approach}
\author{Wen-Long Jin\footnote{Department of Civil and Environmental Engineering, California Institute for Telecommunications and Information Technology, Institute of Transportation Studies, 4000 Anteater Instruction and Research Bldg, University of California, Irvine, CA 92697-3600. Tel: 949-824-1672. Fax: 949-824-8385. Email: wjin@uci.edu. Corresponding author}}

\maketitle

\begin{abstract}
In transportation and other types of facilities, various queues arise when the demands of service are higher than the supplies, and many point and fluid queue models have been proposed to study such queueing systems. However, there has been no unified approach to deriving such models, analyzing their relationships and properties, and extending them for networks. In this paper, we derive point queue models as limits of two link-based queueing model: the link transmission model and a link queue model. With two definitions for demand and supply of a point queue, we present four point queue models, four approximate models, and their discrete versions. We discuss the properties of these models, including equivalence, well-definedness, smoothness, and queue spillback, both analytically and with numerical examples. We then analytically solve Vickrey's point queue model and stationary states in various models. We demonstrate that all existing point and fluid queue models in the literature are special cases of those derived from the link-based queueing models. Such a unified approach leads to systematic methods for studying the queueing process at a point facility and will also be helpful for studies on stochastic queues as well as networks of queues.
 
\end{abstract}


\section{Introduction}
 When the demands of service are higher than the supplies  at such facilities as road networks,  security check points in airports, supply chains, water reservoirs, document processors, task managers, and computer servers, there arise queues of vehicles,  customers, commodities, water, documents, tasks, and programs, respectively. Many strategies have been developed to control, manage, plan, and design such queueing processes so as to improve their performance in safety, efficiency, environmental impacts, and so on. 
 
Since last century, many queueing models have been proposed to understand the characteristics of such  systems, including waiting times, queue lengths, and stationary states and their stability \citep{kleinrock1975queueing,newell1982queueing}. These models can be represented by arrival and service of individual customers or accumulation and dissipation of continuum fluid flows; they can be continuous or discrete in time; and the arrival and service patterns of queueing contents can be random or deterministic. In addition, a queue can be initially empty or not, the storage capacity of a queue can be finite or infinite, there can be single or multiple servers, and the serving rule can be First-In-First-Out or priority-based. 
 Traditionally, queueing theories concern with random arrival and service processes of individual customers \citep{lindley1952theory,asmussen2003applied,kingman2009first}.
 
In recent years, fluid queue models have gained much attention in various disciplines \citep{kulkarni1997fluid}. In these models, the dynamics of a queueing system are described by changes in continuum flows of queueing contents. Such models were first proposed for dam processes in 1950's \citep{moran1956probability,moran1959theory}. Since then, both discrete and continuous versions of fluid queue models have been extensively discussed with random arrival and service patterns.
In the transportation literature, point queue models were introduced to study the congestion effect of a bottleneck with deterministic, dynamic arrival patterns \citep{vickrey1969congestion}. Such models have been applied to study dynamic traffic assignment problems \citep{drissi1992dynamic,kuwahara1997decomposition,li2000reactive}.
In \citep{nie2005comparative}, Vickrey's point queue model was compared with other network loading models.
In \citep{armbruster2006model}, a point queue model for a supply chain was proposed and shown to be consistent with a hyperbolic conservation law.
In \citep{ban2012continuous}, continuous point queue models were presented as differential complementarity systems. In \citep{han2013partial}, Vickrey's point queue model was formulated with the help of the LWR model \citep{lighthill1955lwr,richards1956lwr}. 
In existing studies on point queue models, the downstream service rate is usually assumed to be constant, and the storage capacity of the point queue infinite. Note that fluid queue models can also be considered as point queue models, but with a finite storage capacity and general service rates, since the dimension of a server or a reservoir is also assumed to be infinitesimal. However, there exists no unified approach to deriving such point queue models, analyzing their relationships and properties, and extending them for networks.

In this study we attempt to fill this gap by presenting a unified approach, which is based on two observations. The first observation is that many traffic flow models can be viewed as queueing models, since traffic flow models describe the dynamics associated with the accumulation and dissipation of vehicular queues. For examples, the Cell Transmission Model (CTM) describes dynamics of cell-based queues \citep{daganzo1995ctm,lebacque1996godunov}; the Link Transmission Model (LTM) describes dynamics of link-based queues \citep{yperman2006mcl,yperman2007link}; and the Link Queue Model (LQM) also describes dynamics of link-based queues \citep{jin2012_link}. 
Among these models, continuous versions of LTM and CTM are equivalent to but different formulations of network kinematic wave models, which can admit discontinuous shock wave solutions and can be unstable in diverge-merge networks \citep{jin2014_ltm}, but LQM admits smooth solutions and is stable \citep{jin2012_link}. 
The second observation is that a point facility can be considered as the limit of a link: if we shrink the link length but maintain the storage capacity, then we obtain a point queue.
Therefore we can derive point queue models as limits of the link-based queueing models: LTM and LQM. In this approach, point queue models inherit two critical components from both LTM and LQM: we first define demand and supply of a point queue, and then apply macroscopic junction models, which were originally proposed in CTM and later adopted for both LTM and LQM, to calculate boundary fluxes from upstream demands and downstream supplies. This approach enables us to derive all existing point and fluid queue models and their generalizations, discuss their relationships, and analyze their properties.

The rest of the paper is organized as follows. After reviewing LTM and LQM in Section 2, we take the limits of the demand and supply of a link to obtain two demand formulas and two supply formulas for a point queue in Section 3. Further with different combinations of demand and supply formulas, we derive four continuous point queue models and their discrete versions.  In Section 4, we present four approximate point queue models and their discrete versions. In Section 5, we present some analytical solutions of these models. In Section 6, we present numerical results. Finally in Section 7, we summarize the study and discuss future research topics. 

\section{Review of two link-based queueing models}
In this section, we consider a homogeneous road link, shown in \reff{road_segment}, whose length is $L$ and number of lanes $N$. Vehicles on the road segment form a link queue. We assume that the road has the following triangular fundamental diagram of flow-density relation \citep{munjal1971multilane,haberman1977model,newell1993sim}:
\bqn
q&=&Q(k)=\min\{V k, (NK-k)W\},
\eqn
where $k(x,t)$ is the total density at $x$ and $t$, $q(x,t)$ the total flow-rate, $K$ the jam density per lane, $V$ the free-flow speed, and $-W$ the shock wave speed in congested traffic. We denote the harmonic mean of $V$ and $W$ by $U=\frac{VW}{V+W}$. We also denote $T_1=\frac LV$, which is the free-flow traverse time, $T_2=\frac LW$, which is the shock wave traverse time, and $T_3=\frac LU$. Generally $T_1<T_2<T_3$.
Thus for each lane, the  critical density is $\bar K=\frac{W}{V+W}K=\frac UV K$, and the capacity is $C=V \bar K=U K$. 
We denote the storage capacity of the road link by $\La =NLK$ and the initial number of vehicles by $\La_0 =NLk_0$, where the initial density is constant at $Nk_0$ at all locations on the road. Then the total capacity is  $NC= \frac{\La}{T_3} $.

For a road link shown in \reff{road_segment}, we define the upstream and downstream cumulative flows by $F(t)$ and $G(t)$, respectively, and set $G(0)=0$ and $F(0)=\La_0$. Correspondingly the upstream in- and downstream out-fluxes are denoted by $f(t)$ and $g(t)$. Then we have
\bsq \label{flow-flux}
\bqn
\der{}{t} F(t)&=&f(t),\\
\der{}{t} G(t)&=&g(t).
\eqn
\esq
In addition, we denote the number of vehicles on the road by $\r(t)= \int_0^L k(x,t) dx$, which satisfies
\bqn
\r(t)&=&F(t)-G(t). \label{numveh}
\eqn

\bfg
\bc
\includegraphics[width=4in]{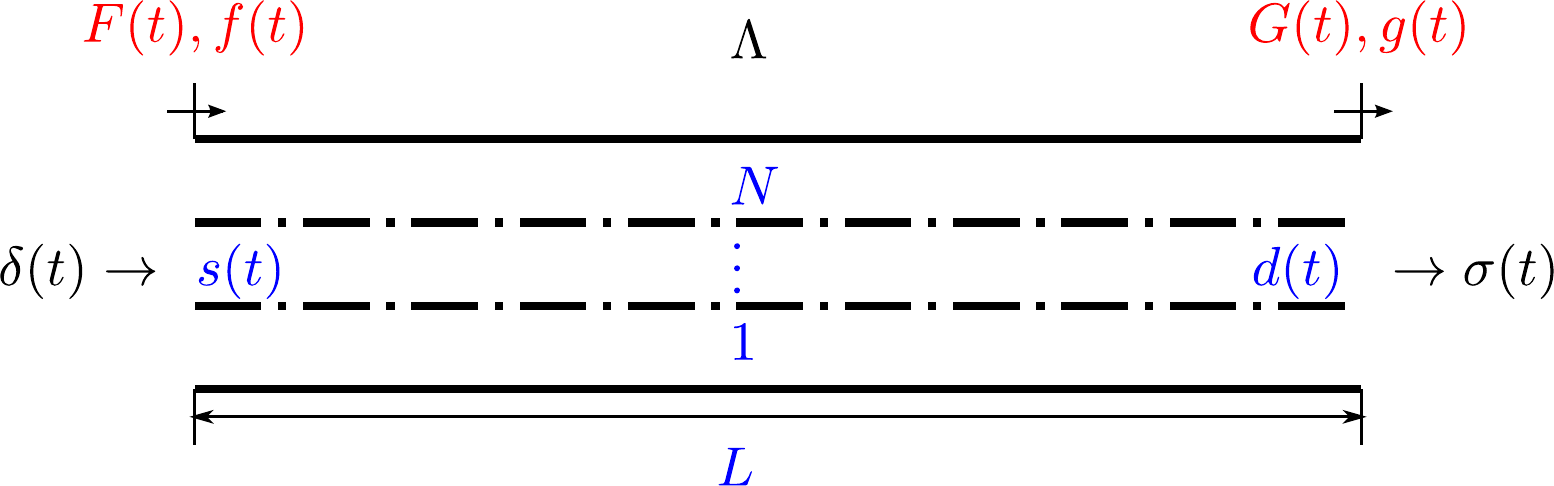}\caption{An illustration of a homogeneous road link}\label{road_segment}
\ec\efg

There are two formulations for a link-based queueing model with two types of state variables: either the accumulation of vehicles, such as $\r(t)$, or the cumulative flows, $F(t)$ and $G(t)$, can be used as state variables. In addition, the boundary fluxes are calculated in two steps: first, a link's demand, $d(t)$, and supply, $s(t)$, are defined from state variables; second,  the flux through a boundary can be calculated from the upstream demand and downstream supply with macroscopic junction models:
 \bsq \label{boundaryflux}
 \bqn
 f(t)&=&\min\{\delta(t),s(t)\},\\
 g(t)&=&\min\{d(t),\sigma(t)\},
 \eqn
 \esq 
where $\delta(t)$ is the origin demand, and $\sigma(t)$ the destination supply.  
LTM and LQM differ in the definitions of demand and supply functions and the corresponding state variables related to vehicle accumulations.

\subsection{The link transmission model (LTM)}
In LTM, two new variables related to vehicle accumulation are the link queue size, $\la(t)$, and the link vacancy size, $\gamma(t)$, which are given by \citep{jin2014_ltm}:
\bsq \label{qv-ltm}
\bqn 
\la(t)&=& \cas{{ll} \frac{\La_0}{T_1} t -G(t), & t\leq T_1 \\F(t-T_1)-G(t), & t>T_1}\\
\gamma(t)&=& \cas{{ll} \frac{\La-\La_0}{T_2} t+\La_0-F(t), & t\leq T_2\\ G(t-T_2)+\La-F(t), &t>T_2}
\eqn
\esq
Since both $F(t)$ and $G(t)$ are continuous, then $\la(t)$ and $\gamma(t)$ are continuous with $\la(0)=0$ and $\gamma(0)=0$.

The link demand and supply in LTM are defined as \citep{jin2014_ltm}:
\bsq \label{ds-ltm}
\bqn
d(t)&=&\cas{{ll} \min\left\{\frac{\La_0}{T_1} +H(\la(t)), \frac{\La}{T_3}\right\}, & t\leq T_1\\
\min\left\{ f(t-T_1) +H(\la(t)), \frac{\La}{T_3} \right\}, & t>T_1}\\
s(t)&=&\cas{{ll} \min\left\{\frac{\La-\La_0}{T_2}+H(\gamma(t)), \frac{\La}{T_3}\right\}, & t\leq T_2\\  \min\left\{ g(t-T_2) +H(\gamma(t)),\frac{\La}{T_3} \right\}, &t>T_2}
\eqn
\esq
where the indicator function $H(y)$ for $y\geq 0$ is defined as
\bqn
H(y)&=&\lim_{\dt\to 0^+}\frac{y}{\dt}=\cas{{ll}0, &y=0\\\infty, &y>0} \label{indicatorfunction}
\eqn

Then we obtain two continuous formulations of LTM:
\ben
\item When the queue and vacancy sizes, $\la(t)$ and $\gamma(t)$, are the state variables, we have  from \refe{qv-ltm}:
\bsq \label{ltm1}
\bqn
\der{}{t} \la(t)&=& \cas{{ll} \frac{\La_0}{T_1}  -g(t), & t\leq T_1 \\ f(t-T_1)-g(t), & t>T_1}\\
\der{}{t}  \gamma(t)&=& \cas{{ll} \frac{\La-\La_0}{T_2} -f(t), & t\leq T_2\\ g(t-T_2)-f(t), &t>T_2}
\eqn
\esq
where $f(t)$ and $g(t)$ are calculated from \refe{boundaryflux} and \refe{ds-ltm}. In this formulation, the rates of change in both $\la(t)$ and $\gamma(t)$ depend on $\lambda(t)$, $\gamma(t)$, $f(t-T_1)$, and $g(t-T_2)$. 
\item When cumulative flows, $F(t)$ and $G(t)$, are the state variables, we have from \refe{flow-flux}:
\bsq \label{ltm2}
\bqn
\der{}{t} F(t)&=&\min\{\delta(t),s(t)\},\\
\der{}{t}G(t)&=&\min\{d(t),\sigma(t)\}.
\eqn
\esq
where $d(t)$ and $s(t)$ are calculated from \refe{qv-ltm} and \refe{ds-ltm}. In this formulation, the rates of change in $F(t)$ and $G(t)$ depend on $F(t)$, $G(t)$, $F(t-T_1)$, $G(t-T_2)$, $f(t-T_1)$, and $g(t-T_2)$. 
\een
In both formulations, $f(t-T_1)$, and $g(t-T_2)$ can be calculated from the corresponding state variables at $t-T_1$ and $t-T_2$. In addition, the right-hand sides of both formulations are discontinuous, and LTM is a system of discontinuous ordinary differential equations with delays. Thus their solutions may not be differentiable or smooth.

\subsection{The link queue model (LQM)}
In LQM, $\r(t)$ can be used as a state variable. The link demand and supply in LQM are defined as \citep{jin2012_link}: 
\bsq\label{ds-lqm}
\bqn
d(t)&=&\min\{ \frac{ \r(t)}{T_1}, \frac{\La}{T_3}\},\\
s(t)&=&\min\{\frac{\Lambda-\r(t)}{T_2}, \frac{\La}{T_3}\}.
\eqn
\esq

From \refe{flow-flux}, \refe{numveh}, \refe{boundaryflux}, and \refe{ds-lqm} we obtain the following two continuous formulations of LQM:
\ben
\item When $\r(t)$ is the state variable, we have: 
\bqn
\der{}{t} \r(t)&=&\min\{\delta(t), \frac{\Lambda-\r(t)}{T_2}, \frac{\La}{T_3}\}- \min\{ \frac{ \r(t)}{T_1}, \frac{\La}{T_3},\sigma(t)\} . \label{lqm1}
\eqn
\item When $F(t)$ and $G(t)$ are the state variables, we have:
\bsq \label{lqm2}
\bqn
\der{}{t}F(t)&=&\min\{\delta(t), \frac{\La -(F(t)-G(t))}{T_2}, \frac{\La}{T_3}\},\\
\der{}{t}G(t)&=&\min\{ \frac{F(t)-G(t)}{T_1}, \frac{\La}{T_3}, \sigma(t) \}.
\eqn
\esq
\een
In both formulations, the right-hand sides are continuous, and LQM are systems of continuous ordinary differential equations without delays. Therefore, solutions of LQM are smooth with continuous derivatives of any order.

\section{Four point queue models}

\bfg
\bc
\includegraphics[width=2.5in]{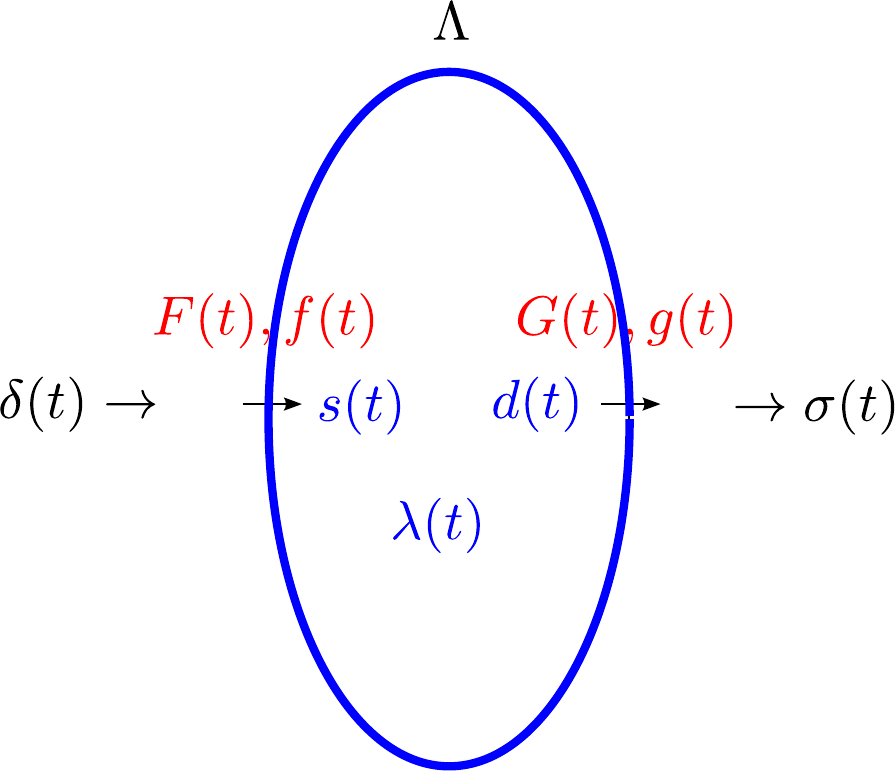}\caption{An illustration of a point queue}\label{point_queue_ill}
\ec\efg

We consider a point facility as a limit of a road link, shown in \reff{point_queue_ill}, with the same storage capacity, $\La$, origin demand, $\delta(t)$, and destination supply, $\sigma(t)$. The macroscopic junction models are also the same as \refe{boundaryflux}. But we let the link length $L\to 0$, and the number of lanes $N=\frac {\La }{LK} \to \infty$.  Parking lots, security check points, water reservoirs, document processors, task managers, and computer servers can all be approximated by point facilities. 
The queue in a point facility is referred to as a point queue.

\subsection{Demand and supply of a point queue}
Based on the observation that a point facility is a limit of a link, we can derive demand and supply of a point queue from \refe{ds-ltm} in LTM and \refe{ds-lqm} in LQM by letting $L\to 0$. 

\begin{lemma} \label{lemma:qv} When $L\to 0$, from the queue and vacancy sizes of LTM, \refe{qv-ltm}, we obtain the queue and vacancy sizes of a point queue as (for $t\geq 0^+$)
\bsq \label{pq-queuesize}
\bqn
\la(t)&=&F(t)-G(t),\\
\gamma(t)&=&\La -\la(t).
\eqn
\esq
\end{lemma}
{\em Proof}. From (\ref{qv-ltm}a) we have $\la(T_1)=\La_0-G(T_1)$ and $\la(t)=F(t-T_1)-G(t)$ for $t>T_1$. When $L\to 0$, we have $T_1\to 0^+$, $\la(0^+)=\La_0 $ and $\la(t)=F(t)-G(t)$ for $t>0$. Thus $\la(t)=F(t)-G(t)$ for $t\geq 0^+$.

Similarly from (\ref{qv-ltm}b) we have $\gamma(T_2)=\La -F(T_2)$ and $\gamma(t)=G(t-T_2)+\La -F(t)$ for $t>T_2$. When $L\to 0$, we have $T_2\to 0^+$, $\gamma(0^+)=\La -\La_0 $ and $\gamma(t)=\La -(F(t)-G(t))$ for $t>0$. Thus $\gamma(t)=\La -\la(t)$ for $t\geq 0^+$.
\eop

From Lemma \ref{lemma:qv}, we can see that, in a point queue, $\la(t)=\r(t)=\La-\gamma(t)$. Thus we can use $\la(t)$, the queue size, to uniquely represent the vehicle accumulation.

\begin{lemma} When $L\to 0$, from the demand and supply in LTM, \refe{ds-ltm}, we obtain the demand and supply  of a point queue as (for $t\geq 0^+$) 
\bsq \label{ds-ltm-pq}
\bqn
d(t)&=&\delta(t)+H(\la(t)),\\
s(t)&=&\sigma(t)+H(\La -\la(t)).
\eqn
\esq
\end{lemma}
{\em Proof}. When $t>T_1$ and $t>T_2$, we have from \refe{boundaryflux} and \refe{ds-ltm}
\bqs
d(t)&=&\min\{\min\{\delta(t-T_1),s(t-T_1)\}+H(\la(t)),\frac{\La}{T_3}\},\\
s(t)&=&\min\{\min\{d(t-T_2),\sigma(t-T_2)\}+H(\gamma(t)),\frac{\La}{T_3}\}.
\eqs 
Let $L\to 0$, then $T_1, T_2, T_3\to 0^+$, $\gamma(t)\to \La -\la(t)$, and \footnote{Without loss of generality we assume that $\delta(t)$, $\sigma(t)$, $d(t)$, and $s(t)$ are all continuous in time.} 
\bqs
d(t)&=&\min\{\delta(t),s(t)\}+H(\la(t)),\\
s(t)&=&\min\{d(t),\sigma(t)\}+H(\La -\la(t)).
\eqs
Therefore, when $\la(t)<\La $ for $t>0$, then $s(t)=\infty$, and $d(t)=\delta(t)+H(\la(t))$;
when $\la(t)>0$ for $t>0$, then $d(t)=\infty$, and $s(t)=\sigma(t)+H(\La -\la(t))$. Thus we obtain \refe{ds-ltm-pq}.
\eop

\begin{lemma} When $L\to 0$, from the demand and supply in LQM, \refe{ds-lqm}, we obtain the demand and supply of a point queue as (for $t\geq 0^+$)
\bsq\label{ds-lqm-pq}
\bqn
d(t)&=&H(\la(t)),\\
s(t)&=&H(\La-\la(t)).
\eqn
\esq
\end{lemma}
{\em Proof}. When $L\to 0$, $T_1, T_2, T_3\to 0$. Thus $d(t)=0$ when $\r(t)=0$ and $\infty$ otherwise. Similarly, $s(t)=0$ when $\r(t)=\La$ and $\infty$ otherwise. Hence \refe{ds-lqm-pq} is correct. \eop 

Clearly the demand and supply functions derived from LTM are different from those from LQM. They can lead to different point queue models.

\subsection{Point queue models}

To model dynamics of a point queue, we can use the queue length, $\la(t)$, or the cumulative flows, $F(t)$ and $G(t)$, as state variables. Once the demand and supply functions are given, from \refe{flow-flux}, \refe{numveh}, \refe{boundaryflux} we can derive two formulations of a point queue model:
\bi
\item [A.] With $\la(t)$ as the state variable, we have
\bqn
\der{}{t} \la(t)&=&\min\{\delta(t), s(t)\}-\min\{d(t),\sigma(t)\}, \label{A-PQM}
\eqn
which is referred to as A-PQM.
\item [B.] With $F(t)$ and $G(t)$ as the state variables, we have
\bsq \label{B-PQM}
\bqn
\der{}{t} F(t)&=&\min\{\delta(t), s(t)\},\\
\der{}{t} G(t)&=&\min\{d(t),\sigma(t)\},
\eqn
\esq
which is referred to as B-PQM.
\ei

\btb\bc
\begin{tabular}{|l||*{2}{c|}} \hline
\backslashbox{$s(t)$}{$d(t)$} &$\delta(t)+H(\la(t))$&$H(\la(t))$\\\hline\hline
$\sigma(t)+H(\La-\la(t))$ & PQM1 & PQM4 \\\hline
$H(\La-\la(t))$ &PQM3& PQM2 \\\hline
\end{tabular}
\caption{Demand and supply functions for four point queue models (PQMs): $\la(t)=F(t)-G(t)$}\label{four_pqm}
\ec\etb

In the preceding subsection, we obtain from LTM and LQM two different definitions for both demand and supply of a point queue, which are functions of $\la(t)$ or $F(t)$ and $G(t)$. Therefore, as shown in \reft{four_pqm}, we then obtain four combinations of demand and supply functions, which lead to four point queue models. Here PQM1 is the limit of LTM, PQM2 the limit of LQM, but PQM3 and PQM4 are mixtures of the limits of both LTM and LQM.

\begin{definition}
A continuous point queue model is well-defined if and only if $\la(t)$ is always between 0 and $\La $ in A-PQM; or   $G(t)\leq F(t) \leq G(t)+\La $ in B-PQM.
\end{definition}

\begin{theorem}
With demand and supply functions in \reft{four_pqm}, all  four point queue models are well-defined. That is, $\la(t)$ is always between 0 and $\La $ in \refe{A-PQM}, and  $G(t)\leq F(t) \leq G(t)+\La $ in \refe{B-PQM}.
\end{theorem}
{\em Proof}. When $\la(t)=0$, $s(t)=\infty$, and $\der{}t\la(t)=\delta(t)-\min\{d(t),\sigma(t)\}\geq 0$, since $d(t)\leq \delta(t)$; when $\la(t)=\La$, $d(t)=\infty$, and $\der{}t\la(t)=\min\{\delta(t), s(t)\}-\sigma(t) \leq 0$, since $s(t)\leq \sigma(t)$. Similarly from  we can see that $G(t)\leq F(t) \leq G(t)+\La $ in \refe{B-PQM}. \eop 

\begin{theorem} \label{thm:equiv} The four point queue models are equivalent.
\end{theorem}
{\em Proof}. When $0<\la(t)<\La$, we have $d(t)=s(t)=\infty$ in all four point queue models, and $\der{}{t}\la(t)=\delta(t)-\sigma(t)$ for all models. When $\la(t)$ increases from 0 or decreases from $\La$, the demand and supply functions and, therefore, rates of change in $\la(t)$ are different for these models. But since the number of such instants is countable, the integrals of different rates of change at such instants are zero. That is, $\la(t)$ is the same for all models, when they have the same origin demand and destination supply. This leads to the equivalence among all models. \eop

If the storage capacity of a point facility is infinite; i.e., if $\La =\infty$, then its supply $s(t)=\infty$ in all four models. In this case, PQM1 and PQM3 are the same, and so are PQM2 and PM4. In particular, A-PQM1 becomes
\bqn
\der{}{t} \la(t)&=&\max\{- H(\la(t)),\delta(t)-\sigma(t)\}=\cas{{ll} \delta(t)-\sigma(t), & \la(t)>0\\ \max\{0,\delta(t)-\sigma(t)\}, & \la(t)=0 } \label{A-PQM1-inf}
\eqn 
which is referred to as A-PQM1i, and B-PQM1 becomes
\bsq  \label{B-PQM1-inf}
\bqn
\der{}{t} F(t)&=&\delta(t)\\
\der{}{t} G(t) &=&\cas{{ll} \min\{\delta(t),\sigma(t)\}, & F(t)=G(t)\\ \sigma(t), & F(t)>G(t) }
\eqn
\esq
which is referred to as B-PQM1i.
Similarly, we can also obtain the special cases of PQM2.

In A-PQM1i, if we allow $\delta(t)$, $\sigma(t)$, and $\la(t)$ to be random variables, then \refe{A-PQM1-inf} is the same as the traditional fluid queue model \citep[][Equation 1]{kulkarni1997fluid} with the difference between origin demand and destination supply, $\delta(t)-\sigma(t)$, as the drift function.
Also in A-PQM1i, if the destination supply, $\sigma(t)$, is constant, then \refe{A-PQM1-inf} is the same as the continuous version of Vickrey's point queue model \citep[][Equations 2.4 and 2.5 with $t_0=0$]{han2013partial}. 
Therefore, PQM1 and PQM3 are generalized versions of existing continuous fluid queue and Vickrey's point queue models, since it also applies to a queue with a finite storage capacity and arbitrary origin demand and destination supply patterns.

\subsection{Discrete versions}

\btb\bc
\begin{tabular}{|l||*{2}{c|}} \hline
\backslashbox{$s(t)\dt$}{$d(t)\dt$} &$\delta(t)\dt+\la(t)$&$\la(t)$\\\hline\hline
$\sigma(t)\dt+\La-\la(t)$ & PQM1-D & PQM4-D \\\hline
$\La-\la(t)$ &PQM3-D& PQM2-D \\\hline
\end{tabular}
\caption{Demand and supply functions for four discrete point queue models: $\la(t)=F(t)-G(t)$}\label{four_pqm-discrete}
\ec\etb

For four continuous point queue models, we discretize them with a time-step size of $\dt$ such that $ H(y)\dt=y$. Then with the discrete demand and supply functions during from $t$ to $t+\dt$, shown in \reft{four_pqm-discrete}, \refe{A-PQM} and \refe{B-PQM} can be numerically solved.
\ben
\item
With $\la(t)$ as the state variable, we have the following four discrete models by replacing $\der{}{t}\la(t)$ by $\frac{\la(t+\dt)-\la(t)}{\dt}$ in \refe{A-PQM}:
\bqn
\la(t+\dt)&=&\la(t)+\min\{\delta(t)\dt,s(t)\dt\}-\min\{d(t)\dt,\sigma(t)\dt\}. \label{A-PQM-D}
\eqn

\item With $F(t)$ and $G(t)$ as the state variables, we have the following four discrete models by replacing $\der{}{t}F(t)$ by $\frac{F(t+\dt)-F(t)}{\dt}$ and $\der{}{t}G(t)$ by $\frac{G(t+\dt)-G(t)}{\dt}$ in \refe{B-PQM}:
\bsq \label{B-PQM-D}
\bqn
F(t+\dt)&=&F(t)+\min\{\delta(t) \dt, s(t)\dt\},\\
G(t+\dt)&=&G(t)+\min\{d(t)\dt,\sigma(t)\dt\}.
\eqn
\esq
\een

In particular, A-PQM3-D can be written as 
\bqn
\la(t+\dt)&=&\min\{\delta(t)\dt+\la(t),\La \}-\min\{\delta(t)\dt+\la(t),\sigma(t)\dt\}, \label{A-PQM1-D}
\eqn
which is the discrete dam process model in \citep[][Chapter 3]{moran1959theory}. In this sense, PQM3 is the continuous version of Moran's model.

\begin{definition} A discrete point queue model, \refe{A-PQM-D}, is well-defined if and only if $\la(t+\dt)$ is between $0$ and $\La$ for $\la(t)$ between $0$ and $\La$; a discrete point queue model, \refe{B-PQM-D}, is well-defined if and only if $G(t+\dt)\leq F(t+\dt)\leq G(t+\dt)+\La$ for $G(t)\leq F(t)\leq F(t)+\La$.
\end{definition}

Then we have the following theorem.

\begin{theorem} 
PQM1-D and PQM2-D are always well-defined with any $\dt$;
 PQM3-D is well-defined if and only if $\dt\leq \frac{\La}{\sigma(t)}$; and PQM4-D is well-defined if and only if $\dt\leq \frac{\La}{\delta(t)}$.
 Thus when
\bqn
\dt \leq \frac {\La}{\max_t\{\delta(t), \sigma(t)\}}, \label{dtcond}
\eqn
all four discrete models are well-defined.
\end{theorem}
{\em Proof}. For A-PQM1-D, we have
\bqs
\la(t+\dt)&=&\min\{\delta(t)\dt+\la(t),\sigma(t)\dt+\La\}-\min\{\delta(t)\dt+\la(t),\sigma(t)\dt\}.
\eqs
Since $\sigma(t)\dt+\La>\sigma(t)\dt$, there are three cases: (i) when $\delta(t)\dt+\la(t)\leq \sigma(t)\dt$, $\la(t+\dt)=0$; (ii) when $\sigma(t)\dt < \delta(t)\dt+\la(t)<\sigma(t)\dt+\La$, $\la(t+\dt)=\delta(t)\dt+\la(t)-\sigma(t)\dt<\La$; (iii) when $\delta(t)\dt+\la(t)\geq\sigma(t)\dt+\La$, $\la(t+\dt)=\La$.
In all three cases, $\la(t+\dt)\in [0,\La]$, and A-PQM1-D is well-defined.

For A-PQM2-D, we have
\bqs
\la(t+\dt)&=&\min\{\delta(t)\dt+\la(t),\La\}-\min\{\la(t),\sigma(t)\dt\}.
\eqs
Since $\la(t)\leq \min\{\delta(t)\dt+\la(t),\La\}\leq \La$, there are two cases: (i) when $\sigma(t)\dt\leq \la(t)$, $\la(t+\dt)=\min\{\delta(t)\dt+\la(t),\La\}-\sigma(t)\dt \in [0,\La]$; (ii) when $\sigma(t)\dt >\la(t)$, $\la(t+\dt)=\min\{\delta(t)\dt+\la(t),\La\}-\la(t) \in [0,\La]$. Thus A-PQM2-D is well-defined.

For A-PQM3-D, we have
\bqs
\la(t+\dt)&=&\min\{\delta(t)\dt+\la(t),\La\}-\min\{\delta(t)\dt+\la(t),\sigma(t)\dt\}.
\eqs
If $\La<\sigma(t)\dt\leq \delta(t)\dt+\la(t)$, then $\la(t+\dt)=\La-\sigma(t)\dt<0$, and the discrete model is not well-defined. However, if $\La\geq\sigma(t)\dt$ or $\dt\leq \frac{\La}{\sigma(t)}$, there are three cases: (i) when $\delta(t)\dt+\la(t)\leq \sigma(t)\dt$, $\la(t+\dt)=0$; (ii) when $\sigma(t)\dt<\delta(t)\dt+\la(t)<\La$, $\la(t+\dt)=\delta(t)\dt+\la(t)-\sigma(t)\dt\in (0,\La)$; (iii) when $\delta(t)\dt+\la(t)\geq \La$,  $\la(t+\dt)=\La-\sigma(t)\dt\in [0,\La]$. Thus A-PQM3-D is well-defined if and only if $\dt\leq \frac{\La}{\sigma(t)}$.

For A-PQM4-D, we have
\bqs
\la(t+\dt)&=&\min\{\delta(t)\dt+\la(t),\sigma(t)\dt+\La\}-\min\{\la(t),\sigma(t)\dt\}.
\eqs
If $\La<\delta(t)\dt$ and $\sigma(t)\dt$ is very large such that $\delta(t)\dt+\la(t)\leq \sigma(t)\dt+\La$ and $\la(t)\leq \sigma(t)\dt$, then $\la(t+\dt)=\delta(t)\dt>\La$, and the discrete model is not well-defined. However, if $\La\geq\delta(t)\dt$ or $\dt\leq \frac{\La}{\delta(t)}$, there are two cases: (i) when $\sigma(t)\leq \la(t)\leq \min\{\delta(t)\dt+\la(t),\sigma(t)\dt+\La\}$, $\la(t+\dt)=\min\{\delta(t)\dt+\la(t),\sigma(t)\dt+\La\}-\la(t)\in [0,\La]$; (ii) when $\sigma(t)> \la(t)$, $\la(t+\dt)=\min\{\delta(t)\dt,\sigma(t)\dt+\La-\la(t)\}\in[0,\La]$. Thus A-PQM4-D is well-defined if and only if $\dt\leq \frac{\La}{\delta(t)}$.

We can similarly prove that the other formulations with $F(t)$ and $G(t)$ as the state variables are also well-defined under the respective conditions. Furthermore, all discrete models are well-defined when \refe{dtcond} is satisfied.  \eop

For a point queue with an infinite storage capacity with $\La =\infty$, all the discrete models are well-defined for any $\dt$, since \refe{dtcond} is always satisfied. In addition, PQM1-D and PQM3-D are the same and can be simplified as
\bqn
\la(t+\dt)&=& \max\{0,\delta(t)\dt+\la(t)-\sigma(t)\dt\}, \label{A-PQM1-D-inf}
\eqn
or
\bsq \label{B-PQM1-D-inf}
\bqn
F(t+\dt)&=&F(t)+\delta(t) \dt,\\
G(t+\dt)&=&\min\{F(t+\dt), \sigma(t) \dt+G(t)\}.
\eqn
\esq
\refe{A-PQM1-D-inf} is the discrete version of Vickrey's point queue model \citep{nie2005comparative}, which is a special case of both PQM1-D and PQM3-D, Moran's discrete dam process model. In addition, \refe{B-PQM1-D-inf} is another formulation of Vickrey's point queue model.
Similarly, PQM2-D and PQM4-D are the same and can be simplified accordingly.

\section{Approximate point queue models}

In the four continuous point queue models derived in the preceding section, the right-hand sides are continuous (assuming that both $\delta(t)$ and $\sigma(t)$ are continuous) but not differentiable due to non-differentiability of the function $H(y)$. In this section we present four approximate models, which admit smooth solutions. 

\subsection{Continuous versions}

\btb\bc
\begin{tabular}{|l||*{2}{c|}} \hline
\backslashbox{$s_\epsilon(t)$}{$d_\epsilon(t)$} &$\delta(t)+\frac{\la_\epsilon(t)}{\epsilon}$&$\frac{\la_\epsilon(t)}{\epsilon}$\\\hline\hline
$\sigma(t)+\frac{\La-\la_\epsilon(t)}{\epsilon}$ & $\epsilon$-PQM1 & $\epsilon$-PQM4 \\\hline
$\frac{\La-\la_\epsilon(t)}{\epsilon}$ &$\epsilon$-PQM3& $\epsilon$-PQM2 \\\hline
\end{tabular}
\caption{Demand and supply functions for four approximate point queue models: $\la_\epsilon(t)=F_\epsilon(t)-G_\epsilon(t)$}\label{four_pqm_appr}
\ec\etb

If we approximate $H(y)\approx\frac{y}{\epsilon}$ for a very small $\epsilon$, then we can obtain differentiable demand and supply functions, $d_\epsilon(t)$ and $s_\epsilon(t)$, shown in \reft{four_pqm_appr}. Denoting the corresponding state variables by $\la_\epsilon(t)$, $F_\epsilon(t)$, and $G_\epsilon(t)$, we obtain the following two formulations of each approximate model:
\bi
\item [A.] With $\la_\epsilon(t)$ as the state variable, we have
\bqn
\der{}{t} \la_\epsilon(t)&=&\min\{\delta(t), s_\epsilon(t)\}-\min\{d_\epsilon(t),\sigma(t)\}, \label{A-PQM-e}
\eqn
which is referred to as $\epsilon$-A-PQM.
\item [B.] With $F_\epsilon(t)$ and $G_\epsilon(t)$ as the state variables, we have
\bsq \label{B-PQM-e}
\bqn
\der{}{t} F_\epsilon(t)&=&\min\{\delta(t), s_\epsilon(t)\},\\
\der{}{t} G_\epsilon(t)&=&\min\{d_\epsilon(t),\sigma(t)\},
\eqn
\esq
which is referred to as $\epsilon$-B-PQM.
\ei

\begin{theorem} \label{thm:approxmodel}
$\epsilon$-PQM1 and $\epsilon$-PQM2 are always well-defined with any $\epsilon$.
 $\epsilon$-PQM3 is well-defined if and only if $\epsilon\leq \frac{\La}{\sigma(t)}$.
$\epsilon$-PQM4 is well-defined if and only if $\epsilon\leq \frac{\La}{\delta(t)}$.
 Thus all approximate models are well-defined, if $\epsilon$ satisfies
\bqn
\epsilon \leq \frac {\La}{\max_t\{\delta(t), \sigma(t)\}}. \label{epsiloncond}
\eqn
\end{theorem}
{\em Proof}. This theorem can be proved by showing that $\der{}{t}\la_\epsilon(t)\geq 0$ when $\la_\epsilon=0$ and $\der{}{t}\la_\epsilon(t)\leq 0$ when $\la_\epsilon=\La$. The proof is straightforward and omitted here. \eop

Since the right-hand sides of these models are differentiable, their solutions are smooth \citep{coddington1972theory}. Clearly when $\epsilon\to 0$, the approximate models converge to the original point queue models.

Furthermore, when $\La=\infty$, \refe{epsiloncond} is satisfied for any $\epsilon$, all the approximate point queue models with an infinite storage capacity are well-defined according to Theorem \ref{thm:approxmodel}. In addition, $\epsilon$-PQM1 is the same as $\epsilon$-PQM3, which can be written as 
\bqn
\der{}t \la_\epsilon(t)&=&\max\{\delta(t)-\sigma(t), -\frac{\la_\epsilon(t)}{\epsilon}\}, \label{pq-1-approx-inf}
\eqn
which is the $\alpha$-model in \citep{ban2012continuous,han2013partial} if we define $\alpha=\frac 1\epsilon$;
$\epsilon$-PQM2 is the same as $\epsilon$-PQM4, which can be written as
\bqn
\der{}{t} \la_\epsilon(t) &=&\delta(t)- \min\{\sigma(t), \frac{\la_\epsilon(t)}\epsilon\},  \label{pq-lqm1-approx-inf}
\eqn
which is the $\epsilon$-model in \citep{armbruster2006autonomous,fugenschuh2008discrete,han2013partial}.

\subsection{Discrete versions}

\btb\bc
\begin{tabular}{|l||*{2}{c|}} \hline
\backslashbox{$s_\epsilon(t)$}{$d_\epsilon(t)$} &$\delta(t)\dt+\frac{\la_\epsilon(t)}{\epsilon} \dt$&$\frac{\la_\epsilon(t)}{\epsilon}\dt$\\\hline\hline
$\sigma(t)\dt+\frac{\La-\la_\epsilon(t)}{\epsilon}\dt$ & $\epsilon$-PQM1-D & $\epsilon$-PQM4-D \\\hline
$\frac{\La-\la_\epsilon(t)}{\epsilon}\dt$ &$\epsilon$-PQM3-D& $\epsilon$-PQM2-D \\\hline
\end{tabular}
\caption{Demand and supply functions for four discrete approximate point queue models: $\la_\epsilon(t)=F_\epsilon(t)-G_\epsilon(t)$}\label{four_pqm_appr-discrete}
\ec\etb

With a time-step size $\dt$, the demand and supply functions for the four discrete approximate point queue models are shown in \reft{four_pqm_appr-discrete}. We can see that when $\dt=\epsilon$, they are the same as the discrete versions of the original point queue models.

\begin{theorem} All the discrete versions of four approximate point queue models are well-defined if and only if the continuous versions are well-defined and $\dt\leq \epsilon$.
\end{theorem}
{\em Proof}. The proof is straightforward and omitted here. \eop

\section{Analytical solutions}
For the four point queue models and their approximate versions, the dynamics are driven by the origin demand, $\delta(t)$, and the destination supply, $\sigma(t)$. In general, as nonlinear ordinary differential equations, they cannot be analytically solved with arbitrary demand and supply patterns. In this section, we first analytically solve Vickrey's point queue model and then solve all point queue models with constant origin demands and destination supplies.

\subsection{Analytical solutions of Vickrey's point queue model}
Even though the destination supply is generally constant in Vickrey's point queue model, we still refer to \refe{A-PQM1-inf} and \refe{B-PQM1-inf} as Vickrey's point queue model with a general destination supply. When $\La_0=0$, its analytical solution is given in the following theorem.

\begin{theorem} \label{thm:vickreypq}
When $\La_0=0$, Vickrey's point queue model, \refe{B-PQM1-inf}, is solved by 
\bsq
\bqn
F(t)&=&\int_0^t \delta(\tau)d\tau,\label{Fsol}\\
G(t)&=&\min_{0\leq \tau \leq t} \{F(\tau)-S(\tau)\}+S(t), \label{Gsol}
\eqn
\esq
where $S(t)=\int_0^T \sigma(\tau)d\tau$,
and \refe{A-PQM1-inf} is solved by
\bqn
\la(t)=F(t)-S(t)-\min_{0\leq \tau \leq t} \{F(\tau)-S(\tau)\}. \label{lasol}
\eqn
\end{theorem}
{\em Proof}. Since $\La_0=0$, $F(0)=G(0)=0$. \refe{Fsol} is straightforward from (\ref{B-PQM1-inf}a). Next we first prove the corresponding discrete version of \refe{Gsol} ($i\geq 0$)
\bqn
G(i \dt)&=&\min_{0\leq j \leq i} \{F(j\dt)-S(j\dt)\}+S(i \dt), \label{Gsol-D}
\eqn
where $S(i\dt)=\sum_{j=0}^{i-1} \sigma(j\dt) \dt$.

When $i=0$, \refe{Gsol-D} is correct. Assuming that \refe{Gsol-D} is true for $i\geq 0$, then from \refe{B-PQM1-D-inf} we have
\bqs
G((i+1)\dt)&=&\min\{F((i+1)\dt),G(i\dt)+\sigma(i\dt) \dt\}.
\eqs
Substituting \refe{Gsol-D} into this equation, we then have
\bqs
G((i+1)\dt)&=&\min\{F((i+1)\dt),\min_{0\leq j \leq i} \{F(j\dt)-S(j\dt)\}+S(i \dt)+\sigma(i\dt) \dt\}\\
&=&\min\{F((i+1)\dt)-S((i+1) \dt),\min_{0\leq j \leq i} \{F(j\dt)-S(j\dt)\}\}+S((i+1) \dt)\\
&=&\min_{0\leq j \leq i+1} \{F(j\dt)-S(j\dt)\}+S((i+1) \dt).
\eqs
Thus \refe{Gsol-D} is correct for $i+1$. From the method of induction, \refe{Gsol-D} is correct for any $i\geq 0$.

For $t\geq 0$, we set $i=\frac t{\dt}$ and $\tau=j\dt$ for $0\leq j\leq i$. Then when $\dt\to 0$, $S(i\dt)\to S(t)$, and \refe{Gsol-D} leads to \refe{Gsol} and \refe{lasol}.
\eop

\begin{corollary} \label{cor:vickreypq}
When the destination supply is constant: $\sigma(t)=\sigma$, then $S(t)=\sigma t$, and the out-flow is 
\bqn
G(t)&=&\min_{0\leq \tau \leq t} \{F(\tau)-\sigma \tau \}+\sigma t, 
\eqn
and the queue length is
\bqn
\la(t)=F(t)-\sigma t-\min_{0\leq \tau \leq t} \{F(\tau)-\sigma \tau\}=\max_{0\leq \tau \leq t} \{F(t)-F(\tau)- (t-\tau)\sigma\}. \label{queuelength-vickrey}
\eqn
\end{corollary}
Note that Corollary \ref{cor:vickreypq} is the same as Theorem 4.10 in \citep{han2013partial}. In addition, \refe{queuelength-vickrey} was also obtained in Lemma 1.2.2 or Corollary 1.5.2 of \citep{leboudec2001network}.

\begin{theorem}
When the destination supply is constant; i.e., when $\sigma(t)=\sigma$, then \citep{li2000reactive},
\bqn
F(t)=G(t+\frac{\la(t)}\sigma). \label{traveltime1}
\eqn
That is, the queueing time for a vehicle entering the queue at $t$ is
\bqn
\pi(t)&=&\frac{\la(t)}\sigma. \label{traveltime}
\eqn
In addition, we have \citep{iryo2007equivalent}
\bqn
(g(t)-\sigma) \cdot \pi(t)&=&0. \label{iryoresult}
\eqn
\end{theorem}
{\em Proof}. When $\la(t)=0$, $F(t)=G(t)$, and \refe{traveltime1} is true. When $\la(t)>0$, for $\tau\in[0,\frac{\la(t)}\tau)$, we have from \refe{B-PQM1-inf} $G(t+\tau)\leq G(t)+\sigma \tau< F(t)\leq F(t+\tau)$, which leads to $\la(t+\tau)>0$ and $g(t+\tau)=\min\{\delta(t+\tau)+H(\la(t+\tau)),\sigma\}=\sigma$. Therefore, $G(t+\frac{\la(t)}\sigma)=G(t)+\la(t)=F(t)$, which is \refe{traveltime1}. Then \refe{traveltime} is true.

Also from $g(t)=\min\{\delta(t)+H(\la(t)),\sigma\}$ we can see that either $\la(t)=0$ and $g(t)\leq \sigma$ or $\la(t)>0$ and $g(t)=\sigma$. Combining \refe{traveltime}, we obtain \refe{iryoresult}.
\eop

\subsection{Stationary states}
In this subsection, we find stationary solutions, $\la(t)=\la$, for the point queue models and their approximations, when both the origin demand and destination supply are constant; i.e., when $\delta(t)=\delta$ and $\sigma(t)=\sigma$. This is the traffic statics problem for a point queue \citep{jin2012statics}. Clearly both $d(t)=d$ and $s(t)=s$ are constant in stationary states. 

\begin{lemma}
In stationary states of a point queue, its in- and out-fluxes are 
\bqn
f=g=\min\{\delta,s\}=\min\{d,\sigma\}=\min\{\delta,s,d,\sigma\}. \label{stationaryflux}
\eqn
\end{lemma}
{\em Proof}. In stationary states, from \refe{boundaryflux} and \refe{A-PQM} we have $f=\min\{\delta,s\}$, $g=\min\{d,\sigma\}$, and $f=g$. Thus we have \refe{stationaryflux}. \eop

Thus we have the following stationary states in the four point queue models.

\begin{theorem}\label{thm:ss1}
All four point queue models have the same stationary states:
\ben
\item When $\delta>\sigma$, $\la=\La$. That is, the point queue is full.
\item When $\delta<\sigma$, $\la=0$. That is, the point queue is empty.
\item When $\delta=\sigma$, $\la\in[0,\La]$. That is, the point queue can be stationary at any state.
\een
\end{theorem}
{\em Proof}.
\ben
\item In PQM1, the stationary states satisfy $\min\{\delta,\sigma+H(\La-\la)\}=\min\{\delta+H(\la),\sigma\}$.
\ben
\item When $\delta>\sigma$, $\la=\La$. That is, the point queue is full.
\item When $\delta<\sigma$, $\la=0$. That is, the point queue is empty.
\item When $\delta=\sigma$, $\la\in[0,\La]$. That is, the point queue can be stationary at any state.
\een
\item In PQM2, the stationary states satisfy $\min\{\delta,H(\La-\la)\}=\min\{H(\la),\sigma\}$.
\ben
\item When $\delta>\sigma$, $\la=\La$ when $\sigma=0$, and there is no solution for $\la$ when $\sigma>0$. However, from its discrete version we have $\la=\La-\sigma \dt \to \La$ when $\dt\to 0$.
\item When $\delta<\sigma$, $\la=0$ when $\delta=0$, and  there is no solution for $\la$ when $\delta>0$. However, from its discrete version we have $\la=\delta \dt \to 0$ when $\dt\to 0$.
\item When $\delta=\sigma$, $\la\in[0,\La]$. 
\een
\item In PQM3, the stationary states satisfy $\min\{\delta,H(\La-\la)\}=\min\{\delta+H(\la),\sigma\}$.
\ben
\item When $\delta>\sigma$, $\la=\La$ when $\sigma=0$, and there is no solution for $\la$ when $\sigma>0$. However, from its discrete version we have $\la=\La-\sigma \dt \to \La$ when $\dt\to 0$.
\item When $\delta<\sigma$, $\la=0$. 
\item When $\delta=\sigma$, $\la\in[0,\La]$. 
\een
\item In PQM4, the stationary states satisfy $\min\{\delta,\sigma+H(\La-\la)\}=\min\{H(\la),\sigma\}$.
\ben
\item When $\delta>\sigma$, $\la=\La$. 
\item When $\delta<\sigma$, $\la=0$ when $\delta=0$, and  there is no solution for $\la$ when $\delta>0$. However, from its discrete version we have $\la=\delta \dt \to 0$ when $\dt\to 0$.
\item When $\delta=\sigma$, $\la\in[0,\La]$.
\een
\een
Thus the theorem is proved.
\eop

Theorem \ref{thm:ss1} confirms that all point queue models are equivalent as shown in Theorem \ref{thm:equiv}.

\begin{theorem}\label{thm:ss2}
All four approximate point queue models, where $\epsilon$ satisfies \refe{epsiloncond}, have the different stationary states:
\ben
\item When $\delta>\sigma$, $\la_\epsilon=\La$ in $\epsilon$-PQM1 and $\epsilon$-PQM4, and $\la_\epsilon=\La-\epsilon \sigma$ in $\epsilon$-PQM2 and $\epsilon$-PQM3.
\item When $\delta<\sigma$, $\la_\epsilon=0$ in $\epsilon$-PQM1 and $\epsilon$-PQM3, and $\la_\epsilon=\epsilon \delta$ in $\epsilon$-PQM2 and $\epsilon$-PQM4.
\item When $\delta=\sigma$, $\la_\epsilon\in[0,\La]$ in $\epsilon$-PQM1, $\la_\epsilon\in[\epsilon \delta,\La-\epsilon \sigma]$ in $\epsilon$-PQM2, $\la_\epsilon\in[0,\La-\epsilon \sigma]$ in $\epsilon$-PQM3, $\la_\epsilon\in[\epsilon \delta,\La]$ in $\epsilon$-PQM4.
\een
\end{theorem}
{\em Proof}. The proof is straightforward and omitted. \eop

\begin{corollary} \label{cor:convergence}
When $\epsilon\to 0$, the stationary states of the approximate point queue models converge to those of the original models.
\end{corollary}
{\em Proof}. This can be easily proved by comparing Theorems \ref{thm:ss1} and \ref{thm:ss2}.\eop

Corollary \ref{cor:convergence} confirms that the approximate point queue models converge to the original ones when $\epsilon\to 0$.

\section{Numerical results}
In this section we simulate the dynamics of a point queue during a period of two hours ($t\in[0,2]$ hr). The origin demand is $\delta(t)=\max\{2000 \sin (\pi t),1000\}$ vph, and the destination supply $\sigma(t)=1200$ vph. Unless otherwise stated, the storage capacity $\La=200$ veh. We assume that the queue is initially empty: $\La_0=0$.

\subsection{Comparison of the four point queue models and their approximations}
In this subsection, we use the discrete models in Section 3.3 to numerically solve the four models with $\dt=0.01$ hr, which satisfies \refe{dtcond}. Both the overall and zoom-in plots of queue lengths in these models are shown in \reff{four_PQM}. From the figure we have the following observations. (i) From \reff{four_PQM}(a), queue lengths are always between 0 and 200 veh, the storage capacity. This confirms that all four models and their discrete versions are well-defined. (ii) From \reff{four_PQM}(c), the maximum queue length  between 0.6 and 0.8 hr equals $\La-\sigma \dt=200-12=188$ veh for PQM2 and PQM3. From \reff{four_PQM}(d), the minimum queue length after 1.8 hr equals $\delta \dt=10$ veh for PQM2 and PQM4. This is consistent with the discussions in the proof of Theorem \ref{thm:ss1}. (iii) The four models provide different numerical solutions: when the queue starts to accumulate but before its size reaches maximum, PQM1 and PQM3 are the same, and PQM2 and PQM4 are the same; when the queue size is maximum and the queue starts to dissipate before it is empty, PQM1 and PQM4 are the same, and PQM2 and PQM3 are the same; after the queue size reaches minimum,   PQM1 and PQM3 are the same, and PQM2 and PQM4 are the same again.

\bfg\bc
\includegraphics[width=6in]{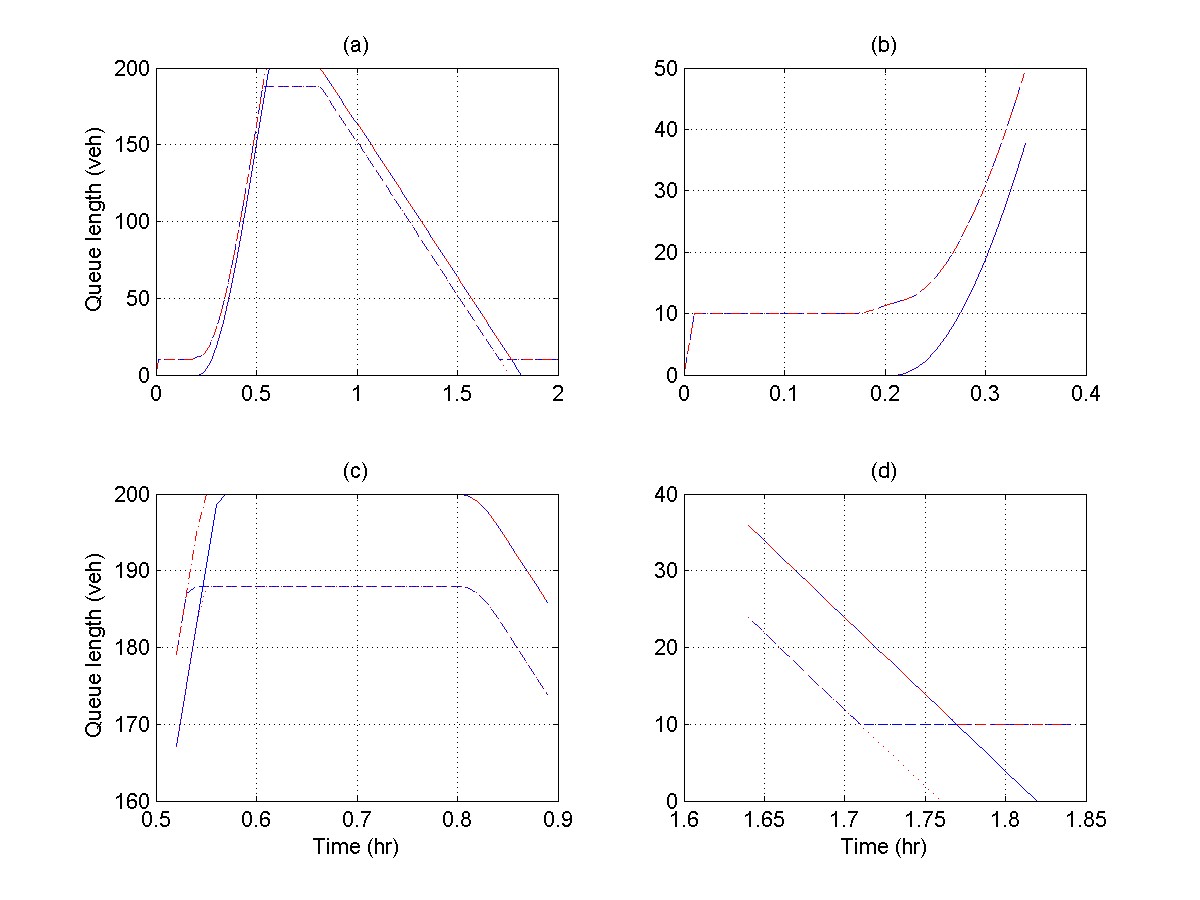}\caption{Comparison of four point queue models: solid lines for PQM1, dashed lines for PQM2, dotted lines for PQM3, and dash-dotted lines for PQM4} \label{four_PQM}
\ec\efg

In \reff{four_PQM_convergence}, we further compare the four models with smaller time-step sizes. Comparing \reff{four_PQM}(a), \reff{four_PQM_convergence}(a), and \reff{four_PQM_convergence}(b), we can see that the four models converge to the same results when $\dt\to 0$. This confirms that all continuous point queue models are equivalent, as stated in Theorem \ref{thm:equiv}. 
The queueing dynamics are as follows: the queue starts to build up at around 0.2 hr, reaches maximum at 0.55 hr, starts to dissipate at 0.8 hr, and disappears at 1.8 hr. Due to the finite storage capacity, the demand that is not satisfied is discarded.
In addition, from \reff{four_PQM_convergence}(b) we can see that the queue length $\la(t)$ is not differentiable when it reaches the maximum and minimum values. Thus, as expected, the solutions of the point queue models are not smooth, even though the origin demand and destination supply are both continuous.

\bfg\bc
\includegraphics[width=6in]{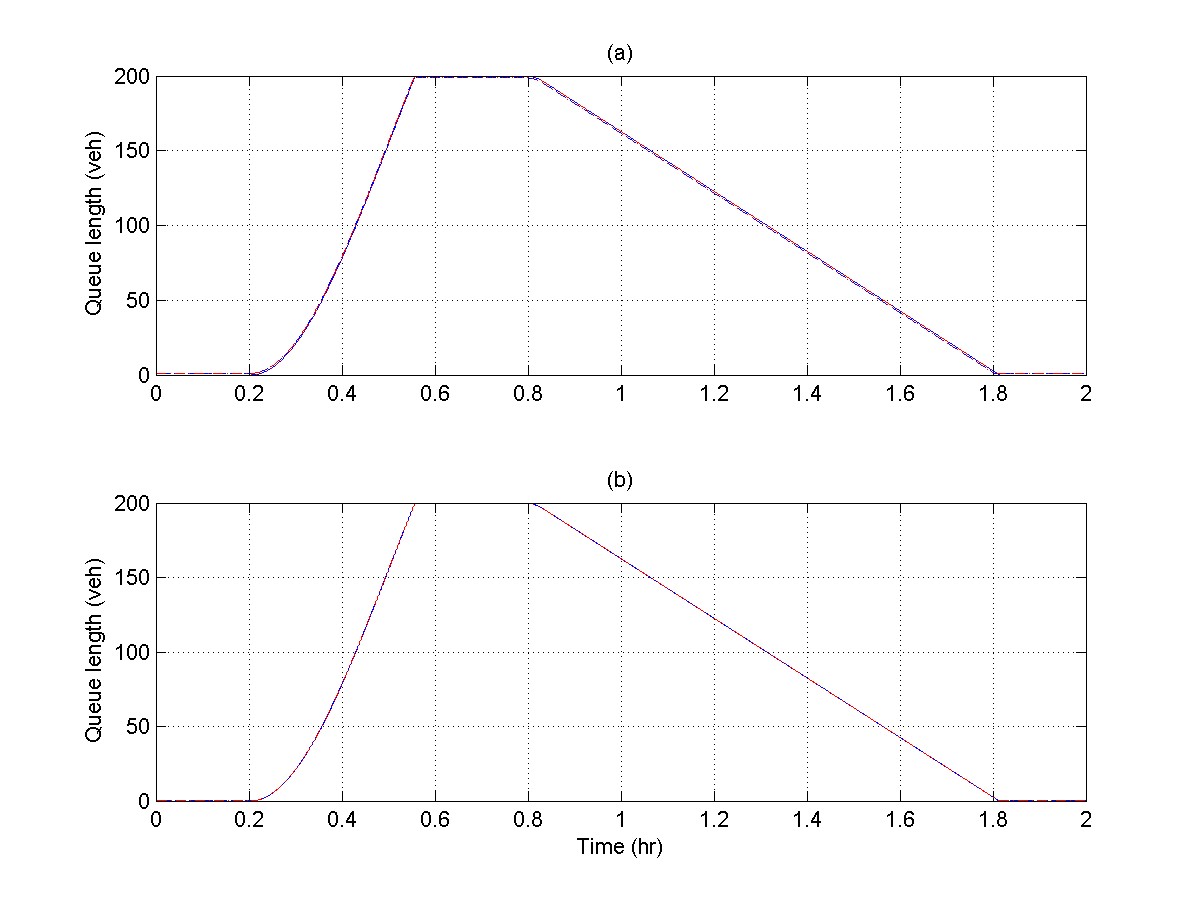}\caption{Comparison of four point queue models: solid lines for PQM1, dashed lines for PQM2, dotted lines for PQM3, and dash-dotted lines for PQM4; (a) $\dt=0.001$ hr; (b) $\dt=0.0001$ hr} \label{four_PQM_convergence}
\ec\efg

Next we solve the discrete versions of four approximate point queue models in Section 4.2 with $\epsilon=0.001$ hr, which satisfies \refe{epsiloncond}. We set $\dt=0.0001\leq \epsilon$ hr. Both the overall and zoom-in plots of the queue length are shown in \reff{four_PQM_approx}, from which we have the following observations: (i) Both continuous and discrete approximate models are well-defined, since the queue length is between 0 and 200 veh; (ii) From \reff{four_PQM_approx}(c), the maximum queue length for $\epsilon$-PQM2 and $\epsilon$-PQM3 is $\La-\epsilon \sigma=198.8$ veh, and that for $\epsilon$-PQM1 and $\epsilon$-PQM4 is 200 veh, as predicted in Theorem \ref{thm:ss2}; (iii) From \reff{four_PQM_approx}(d), the minimum queue length for $\epsilon$-PQM2 and $\epsilon$-PQM4 is $\epsilon \delta=1$ veh, and that for $\epsilon$-PQM1 and $\epsilon$-PQM3 is 0, also as predicted in Theorem \ref{thm:ss2};
(iv) Comparing \reff{four_PQM_approx}(a) and \reff{four_PQM_convergence}(a), we can see that the approximate models converge to the original models when $\epsilon\to0$; (v) From the zoom-in plots we can see that solutions of the approximate point queue models are smooth. 

\bfg\bc
\includegraphics[width=6in]{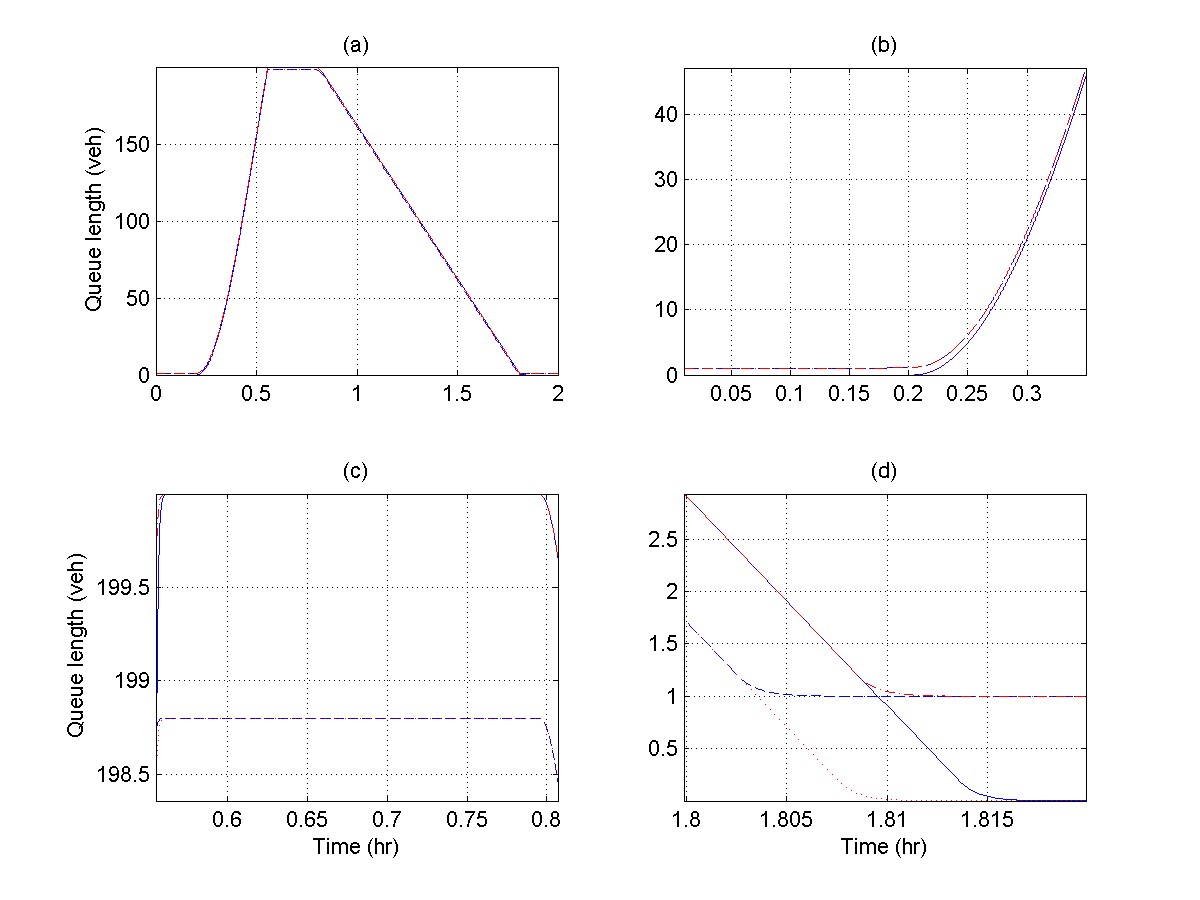}\caption{Comparison of four approximate point queue models: solid lines for $\epsilon$-PQM1, dashed lines for $\epsilon$-PQM2, dotted lines for $\epsilon$-PQM3, and dash-dotted lines for $\epsilon$-PQM4} \label{four_PQM_approx}
\ec\efg

\subsection{Queue spillback effect}
In this subsection we consider a tandem of two point queues: the upstream point queue has an infinite storage capacity, and the downstream one a finite storage capacity $\La_2=200$ veh. Both queues are initially empty. The origin demand and destination supply are the same as in the preceding subsection. 

We denote the queue sizes, demands, and supplies by $\la_i(t)$, $d_i(t)$, and $s_i(t)$ ($i=1,2$), respectively. Then we apply the discrete version of PQM1 in Section 3.3 to simulate the dynamics in such a network:
\bqs
d_1(t)\dt&=&\delta(t)\dt+\la_1(t),\\
s_1(t)\dt&=&\infty,\\
d_2(t)\dt&=&d_1(t)\dt+\la_2(t),\\
s_2(t)\dt&=&\sigma(t)\dt+\La_2-\la_2(t),\\
\la_1(t+\dt)&=&\la_1(t)+\delta(t)\dt-\min\{d_1(t)\dt,s_2(t)\dt\},\\
\la_2(t+\dt)&=&\la_2(t)+\min\{d_1(t)\dt,s_2(t)\dt\}-\min\{d_2(t)\dt,\sigma(t)\dt\}.
\eqs

With $\dt=0.0001$ hr, the simulated queue lengths for both queues are shown in \reff{PQM_D_twoqueues}.  In the tandem of two queues, queue 2 starts to build up at around 0.2 hr; when it reaches its maximum storage capacity at around 0.55 hr, it spills back to queue 1; queue 1 disappears before 1.4 hr, and queue 2 starts to dissipate; and queue 2 does not disappear at 2 hr. 
Therefore the queue spillback effect can be captured in a network of point queue models with a finite storage capacity.

\bfg\bc
\includegraphics[width=4in]{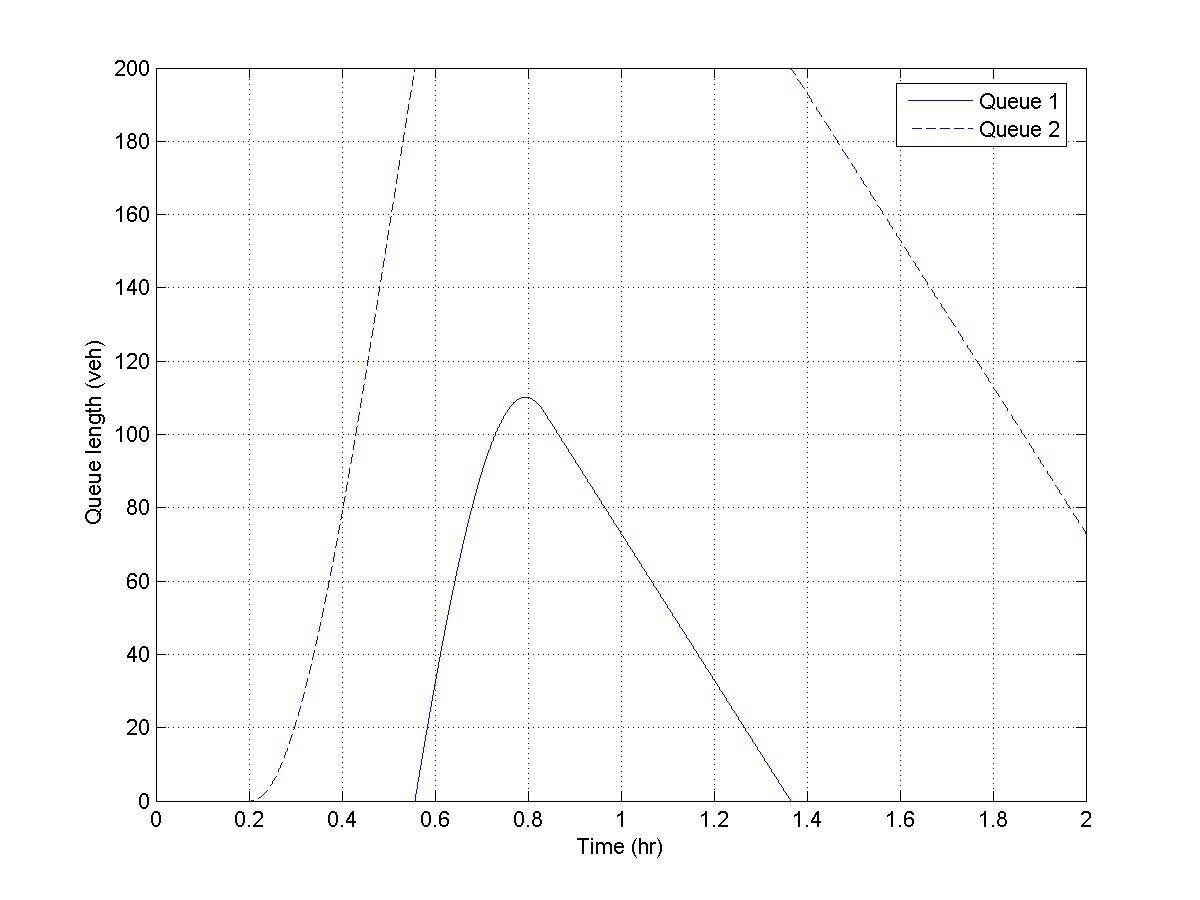}\caption{Queue spillback in a tandem of point queues} \label{PQM_D_twoqueues}
\ec\efg

\section{Conclusions}
In this study we presented a unified approach for point queue models, which can be derived as limits of two link-based queueing model: the link transmission and link queue models. In particular, we provided two definitions of demand and supply of a point queue. The new approach led to the following additional contributions:
\ben
\item From combinations of demand and supply definitions, we presented four point queue models, four approximate models, and their discrete versions. We discussed the following properties of these models: (i) equivalence: the four point queue models are equivalent to each other; (ii) well-definedness: all the continuous and discrete models are well-defined under suitable conditions; and (iii) smoothness: the solutions to the original models are generally not smooth, but those to the approximate models are. These results are obtained theoretically and also verified with numerical simulations. A numerical example also shows that the proposed models with finite storage capacities can capture the queue spillback effect.

\item We analytically solved Vickrey's point queue model, which has an infinite storage capacity and is initially empty. The results are consistent with but more general than those in the literature. We also solved stationary states in the four point models and their approximate versions with constant origin demands and destination supplies. This demonstrates that the new unified formulations enable analytical solutions in special dynamical or stationary cases.

\item In this study we demonstrated that all existing point and fluid queue models in the literature, including Moran's dam process model and Vickrey's point queue model, are special cases of those derived from the link-based queueing models. This further confirms the usefulness of the new unified approach.

\een

Even though in this study point queue models have deterministic origin demands and destination supplies, they can be used to model stochastic queues with random origin demands and destination supplies. For example, for a single queue with an infinite storage capacity, if the origin demand is random and follows a Poisson process, and the destination supply is constant, then this becomes an M/D/1 queue, and some of the results in Section 5.1 still apply. In the future we will be interested in studying stochastic queues with different models derived in this study.

With demand and supply variables, point queue models have the same structure as link-based queueing models. Therefore, they can be integrated into a network of both point and link queues with appropriate macroscopic junction models, which calculate boundary fluxes from upstream demands and downstream supplies \citep{jin2012network}. In the future we will be interested in studying stationary patterns and dynamics in networks of deterministic or stochastic queues, which can arise in computer, air traffic, and other systems. Furthermore, we will also be interested in control, management, scheduling, planning, and design of such queueing networks.


\end{document}